\documentstyle[12pt]{ article}
\begin{document}

\newtheorem{theorem}{Theorem}[section]
\newtheorem{prop}[theorem]{Theorem}
\newenvironment{mtheorem}{\begin{prop}\rm}{\end{prop}}
\newtheorem{coroll}[theorem]{Corollary}
\newenvironment{mcorollary}{\begin{coroll}\rm}{\end{coroll}}
\newtheorem{lemm}[theorem]{Lemma}
\newenvironment{mlemma}{\begin{lemm}\rm}{\end{lemm}}
\newtheorem{anoprop}[theorem]{Proposition}
\newenvironment{mproposition}{\begin{anoprop}\rm}{\end{anoprop}}
\newenvironment{mproof}{\begin{trivlist}\item[]{\em
Proof: }}{\hfill$\Box$\end{trivlist}}
\newenvironment{mdefinition}{\begin{trivlist}\item[]
{\em Definition.}}{\hfill$\Box$\end{trivlist}}
\newtheorem{eg}{\rm\sl \uppercase{Example}}[section]
\newenvironment{example}{\begin{eg}\rm}{\hfill$\Box$\end{eg}}

\newcommand{\ca}{{\cal A}}
\newcommand{\cb}{{\cal B}}
\newcommand{\cc}{{\cal C}}
\newcommand{\cd}{{\cal D}}
\newcommand{\cf}{{\cal F}}
\newcommand{\cg}{{\cal G}}
\newcommand{\ch}{{\cal H}}
\newcommand{\ck}{{\cal K}}
\newcommand{\cl}{{\cal L}}
\newcommand{\cn}{{\cal N}}
\newcommand{\cm}{{\cal M}}
\newcommand{\cs}{{\cal S}}
\newcommand{\ct}{{\cal T}}
\newcommand{\cq}{{\cal Q}}

\newcommand{\csa}{{$C^*$-algebra}}
\newcommand{\nsa}{{non-self-adjoint}}
\newcommand{\cssa}{{$C^*$-subalgebra}}

\def \IR{\hbox{{\rm I}\kern-.2em\hbox{{\rm R}}}}
\def \iR{\hbox{{\sevenrm I\kern-.2em\hbox{\sevenrm R}}}}
\def \IN{\hbox{{\rm I}\kern-.2em\hbox{\rm N}}}
\def \IC{\hbox{{\rm I}\kern-.6em\hbox{\bf C}}}
\def \IQ{\hbox{{\rm I}\kern-.6em\hbox{\bf Q}}}
\def \ZZ{\hbox{{\rm Z}\kern-.4em\hbox{\rm Z}}}

\newcommand{\ga}{{\alpha}}
\newcommand{\gb}{{\beta}}
\newcommand{\gd}{{\delta}}
\newcommand{\gA}{{\Alpha}}
\newcommand{\gB}{{\Beta}}
\newcommand{\gG}{{\Gamma}}
\newcommand{\gD}{{\Delta}}
\newcommand{\gs}{{\sigma}}
\newcommand{\gS}{{\Sigma}}

\newcommand{\ot}{\otimes}
\newcommand{\op}{\oplus}
\newcommand{\pr}{\prime}

\newcommand{\roa}{reflexive operator algebra }

\begin{titlepage}
\begin{center}
{\Large{\bf Homology for Operator Algebras I :\\

Spectral homology for reflexive algebras \\
\vspace{1in}}}

{\Large S.C. Power$^{*}$}
\end{center}
\vspace{1cm}
\begin{center}
         Department of Mathematics\\
Lancaster University\\
England

\end{center}

\vspace{1cm}

\begin{abstract}
A stable homology theory is defined for completely distributive CSL algebras
in terms of the point-neighbourhood homology of the partially ordered set of
meet-irreducible elements of the invariant projection lattice. This specialises
to the simplicial homology of the underlying simplicial complex in the case of
a digraph algebra. These groups are computable and useful. In particular
it is shown
that if the first spectral homology group is trivial then
Schur automorphisms are automatically quasispatial. This motivates
the introduction
of essential Hochschild cohomology which we define by using the point weak star
closure of coboundaries in place of the usual coboundaries.

\end{abstract}
\vspace{.3in}

\footnoterule
\noindent {\footnotesize * Supported by a NATO collaborative research grant. }

\end{titlepage}

\section{Introduction}

A leading theme  in operator algebra is the analysis of
automorphisms and derivations, as well as related Hochschild cohomology.
In the present paper we are concerned largely with a certain tractable
and well-known class of reflexive operator algebras on separable Hilbert
space, namely those with completely distributive commutative invariant
projection lattice. Starting with the partially ordered set of
meet-irreducible projections in the lattice - the crucial ingredient for
the intrinsic spectral representation theory of \cite{orr-scp} and
\cite{scp-tensor}
- we define
 new homology groups, \ $\ H^{sp}_* ({\cal A}),\ $\  which
we call the (integral) spectral
homology of \ $\ \cal A.\ $ In contrast to Hochschild cohomology for Banach
algebras, the spectral homology is often instantly
computable. Furthermore we can obtain a Kunneth formula for
the spectral homology of
spatial tensor products, together with a natural suspension formula.
These formulae follow fairly directly
from the corresponding ones in simplical
homology.

What is less clear, except in finite-dimensions (cf.
Proposition 3.7 below), is the
relationship between  \ $\ H^{sp}_* ({\cal A}),\ $\
and derivations, automorphisms and Hochschild cohomology. However the
main point and result of this paper is that there is a useful connection.
Specifically we show that if \ $\ H^{sp}_1 ({\cal A})\ $\  is trivial then
every
Schur automorphism \ $\ \alpha\ $\
relative to a fixed masa of \ $\ {\cal A}\ $\  is
quasispatial in the sense of Gilfeather and Moore \cite{gil-moo}.
Equivalently, \ $\ \alpha\ $\  lies in the
closure of the inner Schur automorphisms with respect to the point weak
star topology.  Alternatively put, we have the implication

\[
H^{sp}_1 ({\cal A}) = 0 \ \Rightarrow \hbox{ Hoch}_{ess}^1 ({\cal A}) =0
\]

\noindent where \ $\ \hbox{ Hoch}_{ess}^{1} (\cal A)\ $\
is what we call the first
essential Hochschild cohomology group.
The essential Hochschild cohomology arises by replacing the usual space of
coboundaries by their point weak star closure.
The implication  gives an
appealing intuitive coordinate-theoretic understanding of why specific
algebras may have trivial first Hochschild cohomology  (essential or
ordinary, because the two often coincide).
Such examples include nest algebras (cf. Lance \cite{lan}),
and tree
algebras \cite{dpp}, as well as various algebras arising from tensor
products, and the cone and join constructions
considered by Gilfeather and Smith \cite{gil-smi-1},\cite{gil-smi-2}.
We expect that the implication is also valid for the higher order
groups.

Essential Hochschild cohomology does seem to be, in many respects, the most
natural form of Hochschild cohomology for CSL algebras. In the case of
the tridiagonal algebras of the form

\[
{\cal A} = \left[ \begin{array}{cccccc} {\cal L} ({\cal H}_1) & {\cal L} ({\cal
H}_2,\ch_1) & 0 & &\\
0 & {\cal L} ({\cal H}_2) & 0 & &\\
0 & {\cal L} ({\cal H}_2,\ch_3) & {\cal L} ({\cal H}_3) & & \\
& & * & & \\
& & * & * & *\\
& & & & & \ddots \end{array} \right]
\]

\noindent Gilfeather, Hopenwasser and Larson \cite{gil-hop-lar}
have computed the nonzero group
Hoch$^1 ({\cal A}).\ $\  On the
other hand all automorphisms of \ $\ \cal A\ $\  are
quasispatial (see \cite{gil-moo}). For this algebra we have
\ $\ H^{sp}_1 ({\cal A}) =0,\ $\
a fact which is immediately apparent since
this homology group coincides with the simplical
homology of the infinite complex

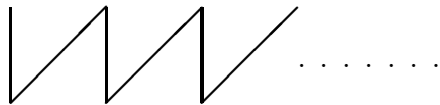
\begin{center}
\setlength{\unitlength}{0.0125in}%
\begin{picture}(120,40)(80,760)
\thicklines
\put( 80,800){\line( 0,-1){ 40}}
\put( 80,760){\line( 1, 1){ 40}}
\put(120,800){\line( 0,-1){ 40}}
\put(120,760){\line( 1, 1){ 40}}
\put(160,800){\line( 0,-1){ 40}}
\put(160,760){\line( 1, 1){ 40}}
\put(200,775){\makebox(0,0)[lb]
{\raisebox{0pt}[0pt][0pt]{\twlrm .  .  .  .  .  .  .  }}}
\end{picture}
\end{center}

\noindent From our perspective it is this vanishing of the first
spectral homology group that is the
underlying reason for the quasispatiality of automorphisms
in this case.

An important ingredient of the proof of the main result is the
representation of maps \ $\ \Phi : \cl(\ch) \to \cl(\ch)\ $\
which are bimodule maps
relative to a maximal abelian self-adjoint subalgebra.
Such maps are sums of elementary ones of the form
\ $\ A \to CAD\ $\ with \ $\ C $\ and \ $D$\ in the masa.
This type of representation goes back to Haagerup \cite{haa} and a convenient
analysis of these issues is given by Smith in \cite{smi}. Although
such representations
play a crucial role in the completely bounded cohomology of
Christensen, Effros and Sinclair \cite{chr-eff-sin}, and in the cohomology
calculations of Gilfeather and Smith, they play a different role in the
proof of our main result.
Using this representation,
together with a local multicativity property for a Schur automorphism
\ $\ \ga \ $\ , we obtain an elementary local representation \ $\ \ga (A) =
CAD\ $\  for all operators
\ $\ \ca\ $\  in the local space \ $\ Q_1 {\cal L} ({\cal H}) Q_3\ $\
associated with a
triple \ $\ Q_1 \prec Q_2 \prec Q_3\ $\  of comparable intervals
of \ \ Lat$\ \ca$. Here \ $\ C\ $\  and \ $\ D\ $\
are determined up to a nonzero scalar. The hypothesis of
vanishing spectral homology is just the condition needed in order
to extend the local
representations to global ones on large enough subalgebras of \ $\ \cal A\ $\
to guarantee quasispatiality. To perform this last step it is necessary
to make use of  the
intrinsic spectral representation of \ \ Lat$\cal A\ $\  obtained
in \cite{orr-scp} together with the additional propositions
obtained below in section 2.  The
proof of the intrinsic representation theorem is quite lengthy and
technical,
and this combines to make the proof of the main result here a long one.
Certainly it
would be interesting to obtain a shorter proof.
Nevertheless, this technical complexity should
not obscure the usefulness of the ideas of the main
result in practice. Indeed, in most concrete examples the
spectral representation and the associated technicalities given in
section 2 are completely apparent. See, for example, the
algebras of section 3.10 below.

The domain of $H^{sp}_{*}$ is the class of reflexive operator
algebras with completely distributive commutative invariant
projection lattice. We leave open the problem of how to extend the domain
to larger classes of \nsa  \ operator algebras. The most obvious task here
is to define spectral homology for general CSL algebras, that is, to
drop the complete distributivity requirement, but we also mention
other possibilities.

In \cite{scp-stable} we develop a quite different homology theory for general
\nsa \ operator algebras. This
is based upon the homology of complexes arising from
partial isometries in the
associated stable algebra and is somewhat $K$-theoretic in nature.
In the final section we comment on
the fact that the spectral homology groups need not coincide with
these partial isometry
homology groups.

The notation and terminology used below is fairly standard
and in the main is coherent with that given in the books of
Davidson \cite{dav-book} and Power \cite{scp-book}.
In particular, on occasion we write $\ \ A(G)\ \ $ for the digraph algebra
(finite-dimensional CSL algebra) associated with the finite directed
transitive graph $\ \ G.\ \ $

The results in this  paper have been incubating and
hibernating over a good few years. I would
like to thank Vern Paulsen and Ken Davidson
for their support in this project,
and also David Larson and Roger Smith for excellent hospitality
in January 1993 when this work was all but completed.

\section{Spectral invariants }

All operator algebras discussed below are assumed
to act on a separable Hilbert space.
A CSL algebra is a \roa  $\ca$ for which the lattice of invariant
(self-adjoint) projections $\cl = $\ Lat$\ca$  is a set of commuting
projections. A CDCSL algebra is a CSL algebra for which $\cl$ is
completely distributive as an abstract lattice. In many respects
this subclass is the most tractable family of \nsa  \ operator algebras.
See, for example, the discussions in
\cite{arv}, \cite{lam}, \cite{dav-book}
\cite{scp-tensor}, \cite{hop},
\cite{dav-scp}, \cite{dpp}, \cite{dav-interpolating} and
\cite{dav-pit}.

The spectral representation theorem for (lattices of)  CDCSL algebras,
which is given in
Power \cite{scp-tensor} and Orr and Power \cite{orr-scp}, uses
as underlying partially
ordered space the set $M(\cl)$ of meet-irreducible elements
of $\cl$, excluding the identity operator. The partial ordering
is the natural one : $L_1 \le L_2$ if and only if $L_1L_2 = L_1$.
By way of illustration consider the 4-cycle digraph algebra \ \ $\ca$\
in $M_4(\IC)$ which is associated with the pattern
\[
\left(
\begin{array}{cccc}
* & 0 & * & *\\
0 & * & * & *\\
0 & 0 & * & 0\\
0 & 0 & 0 & *
\end{array}
\right).
\]
In particular $\ca$ contains the 4 minimal projections
$e_{11}, e_{22}, e_{33}, e_{44}$.
Consider the meet-irreducible projections $E_{ii} = \sup \{L: L \in
\hbox{ Lat}
{\cal A}, \ L e_{ii} =0 \}$.
Thus \ $\ E_{11} = e_{22}, \ E_{22} = e_{11}, \ E_{33} =
e_{11} + e_{22} + e_{44}, \ E_{44} = e_{11} + e_{22} + e_{33}\ $\ . In fact
\ $\ M({\cal L}) = \{E_{ii} : 1 \leq i \leq 4\}\ $\  and the partially ordered
set \ $\ (M ({\cal L}), \leq)\ $\  can be viewed as a copy of the 4-cycle
digraph

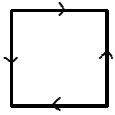
\begin{center}
\setlength{\unitlength}{0.0125in}%
\begin{picture}(45,45)(120,715)
\thicklines
\multiput(143,720)(-0.42857,-0.42857){8}{\makebox(0.4444,0.6667){\sevrm .}}
\multiput(140,717)(0.50000,-0.33333){7}{\makebox(0.4444,0.6667){\sevrm .}}
\multiput(120,737)(0.50000,-0.33333){7}{\makebox(0.4444,0.6667){\sevrm .}}
\multiput(123,735)(0.40000,0.40000){6}{\makebox(0.4444,0.6667){\sevrm .}}
\multiput(165,737)(-0.33333,0.50000){7}{\makebox(0.4444,0.6667){\sevrm .}}
\multiput(163,740)(-0.42857,-0.42857){8}{\makebox(0.4444,0.6667){\sevrm .}}
\put(123,757){\line( 1, 0){ 40}}
\put(163,757){\line( 0,-1){ 40}}
\put(163,717){\line(-1, 0){ 40}}
\put(123,717){\line( 0, 1){ 40}}
\multiput(143,755)(0.40000,0.40000){6}{\makebox(0.4444,0.6667){\sevrm .}}
\multiput(145,757)(-0.33333,0.50000){7}{\makebox(0.4444,0.6667){\sevrm .}}
\end{picture}
\end{center}

\noindent This digraph (the convention is to omit the edges with
a single vertex)
is  more usually associated with $\ca$
in terms of the minimal interval projections of \ $\cl$\
(the atoms of \ $\cl$\ ) with the algebraic partial ordering :
$e_{ii} \prec e_{jj}$ if and only if $e_{ii}\ca e_{jj} = e_{ii}M_4e_{jj}$.
However, for algebras that are not purely atomic,
such as the spatial tensor product $\ca \otimes \ct $, where
$\ct$ is a continuous nest algebra, this kind of simple
association is not available.

For many specific algebras, such as those given at the end of section 3,
it is usually straightforward to
identify the invariant
$M(\cl)$ and its partial ordering. The following proposition is also useful.
\vspace{.3in}

\noindent {\bf Proposition 2.1}\it \  Let Alg({\cal L}) be
a completely distributive CSL algebra. Then

\noindent (i) $M({\cal L})$ is totally ordered if and only if
Alg${\cl}$ is a nest algebra.\\
(ii) $M({\cal L})$ is path connected, when viewed as an infinite graph,
if and only if the commutant of $\ca$ is trivial.\\
(iii) If $\cl_1 \otimes \cl_2$ is the spatial tensor product
then M($\cl_1 \otimes \cl_2)$
= M($\cl_1) \times$  M($\cl_2$) with the product order.\rm

\begin{mproof} (i) This follows from Theorem  2.2 below.

(ii) See Proposition 5.1 of \cite{scp-tensor}.

(iii) The corresponding product structure for nonzero join-irreducible
 elements in a spatial tensor product
has already been given in Lemmas 4.3 and 4.4 of \cite{scp-tensor}. A proof
 of (iii) follows on applying this to the complementary
lattice of $\cl_1 \ot \cl_2$. For completeness we give an
independent proof below in
the proof of the Kunneth formula.
\end{mproof}

The following construction of a CSL algebra in terms of a quadruple - the
order-measure-multiplicity type -  will feature in the fundamental
representation theorem below.

Let $X$ be a partially ordered set with the order topology
generated by the semi-open
order intervals $[x,y)$, let $\mu$ be a Borel measure on $X$,
let $\{E_k\}$ be a
measurable partition of $X$ (by non-null sets),
and let $(n_k)$ be a sequence of
distinct numbers in $\IN \cup \{\infty\}$. The range
of the index $k$ may be finite.
To the quadruple
${\cal X} = (X, \mu, \{E_k\}, (n_k))$ associate the CSL
projection lattice
${\cal L}({\cal X})$ on the Hilbert space
\[
{\cal H} = \sum_k \oplus (L^2(E_k,\mu)\otimes\IC^{n_k})
\]
consisting of the projections
\[
\tilde{P}(B) = \sum_k \oplus (P(B) \otimes I_{n_k})|_{L^2(E_k,\mu) \otimes \;
\IC^{n_k}},
\]
where $B$ is a decreasing Borel subset of $X$, with
associated projection $P(B)$
on $L^2(X,\mu)$.
Of course the resulting lattice need not be completely distributive.
\vspace{.3in}

\noindent {\bf Theorem 2.2 }\  \cite{orr-scp}\em \ Let ${\cal L}$
be a completely distributive
separably acting commutative projection lattice. Let $X$ be the partially
ordered set ${\cal M}({\cal L})$ of meet-irreducible projections of ${\cal L}$,
excepting the identity operator, carrying the order topology generated by the
semi-intervals $[E,F)$. Then there is a
Borel probability measure $\mu$ on $X$, and a
multiplicity partition $\{E_k\},
(n_k),$ such that ${\cal L}$ is unitarily equivalent to the lattice ${\cal
L}({\cal X})$ associated with the quadruple ${\cal X} = (X,\mu, \{E_k\},
(n_k))$.\rm
\vspace{.3in}

The topology on \ $\ M({\cal L})\ $\  plays a minor role in the theorem above,
since it serves merely to locate the sigma algebra on which \ $\ \mu\ $\  is
defined. However, in what follows a topology on \ $\ M({\cal L})\ $\  will be
used in the definition of the spectral homology in terms of a
point-neighbourhood homology and so we now give some attention to this.

The default topology on \ $\ M({\cal L})\ $\  is defined to be that which is
generated by the semi-intervals \ $\ [L, I)\ $\  and their complements \ $\ [L,
I)^c\ $\ , where \ $\ L\ $\  ranges over all projections in the
lattice \ $\ \cal L.\ $\ In
particular, we can take a neighbourhood base for a projection \ $\ M\ $\  in
\ $\ M({\cal L})\ $\  to be the family of
subsets of \ $\ M({\cal L})\ $\  of the form

\[
[M, I) \cap [L_1, I)^c \cap \dots \cap [L_n, I)^c
\]

\noindent where \ $\ L_1, \dots, L_n\ $\
belong to \ $\ \cal L\ $\  and \ $\ n $\
is arbitrary. This
topology is an algebra of sets, and in \cite{orr-scp} it was shown that for a
given separating vector \ $\ e\ $\  for \ $\ \cal L\ $\
there is a unique finitely
additive measure \ $\ \mu\ $\  on this algebra
such that for all \ $\ L_1, \dots,
L_n\ $\  in \ $\ \cal L\ $\  we have

\[
\mu ([L_1, I)^c \cap \dots \cap [L_n, I)^c)\ \  =\ \  <L_1 L_2 \dots L_n e, e>.
\]
\vspace{.3in}

\noindent {\bf Proposition 2.3}\it \  \ If \ $\ B \subseteq
M ({\cal L})\ $\  is open then \ $\ \mu (B)
>0\ $\ .\rm

\begin{mproof} \ Let \ $\
B = [M, I) \cap A\ $\  where \ $\ A = [L_1, I)^c \cap \dots
\cap [L_n, I)^c$. Then, for each \ $\ k,\ $\  \
$\ M^{\bot} L_k \ne 0\ $\ ,
and so \ $\ M_k =
M + M^{\bot}L_k \ $\  is strictly larger than \ $\ M$. On the other hand if
\ $\ \mu (B) = 0\ $\  then \ $\ \mu ([M, I)^c \cap A) = \mu (A)\ $\
and so from the
definition of \ $\ \mu\ $\  and the fact that \ $\ e\ $\  is
separating we have \ $\ M L_1
L_2 \dots L_n = L_1 L_2 \dots L_n$. But now \ $\ M = M_1 M_2 \dots M_k\ $\
contrary to the meet-irreducibility of \ $\ M$.
\end{mproof}
\vspace{.3in}

\noindent {\bf Proposition 2.4}\it \ \  Let \ $\ A, B\ $\  be Borel
subsets of \ $\ M({\cal L})\ $\ , with
projections \ $\ P(A), P(B)\ $\  in \ $\ \cal L$. Then \ $\ \tilde{P}(A)
\prec \tilde{P}(B)\ $\
if and only if the set \ $\ \{(x,y): x \in A, y \in B, x \not\le y\}\ $\  has
product measure zero.\rm

\begin{mproof} Immediate from the definition of \ $\ \tilde{P}(A)
\ $and $\ \tilde{P}(B)$
\end{mproof}

Let us
write \ $\ A \prec B\ $\  if \ $\ x \le y\ $\
for all \ $\ x\ $\  in \ $\ A\ $\
and \ $\ y\ $\  in \ $\ B\ $\ .
For the basic sets of the topology on $\ M(\cl)$ \ we have the
following improvement of the last
proposition.
\vspace{.3in}

\noindent {\bf Proposition 2.5}\it\  \ Let \ $\ A = [M, I) \cap C\ $\
where \ $\ C = [L_1, I)^c \cap
\dots \cap [L_n, I)^c\ $\  and
let \ $\ B = [M^{\prime}, I) \cap C^\prime\ $\  with
\ $\ C^\prime = [L_1^\prime, I)^c \cap \dots \cap [L_n^\prime, I)^c.\ $\  Then
\ $\ \tilde{P}(A) \prec \tilde{P}(B)\ $\  if and only if \ $\ A \prec B$.\rm

\begin{mproof}  That \ $\ A \prec B\ $\  implies \ $\ \tilde{P}(A) \prec
\tilde{P}(B)\ $\  is
routine. For the converse direction notice that it is sufficient to
prove that \ $\ M \leq M^\prime$. Indeed, if this is so then we deduce that
\ $\ M_1 \leq M^\prime\ $\  for all \ $\ M_1\ $\  in \ $\ A\ $\  by
noting that \ $\ [M_1, I) \cap C
\prec B.\ $\  But if \ $\ M \not\le M^\prime\ $\
then \ $\ M^\prime \in [M, I)^c\ $\  and

\[
\tilde{P} ([M, I) \cap C)) \prec \tilde{P} ([M^\prime, I) \cap C^\prime \cap
[M, I)^c)
\]

\noindent where both these projections are nonzero, by Proposition 2.3. By
Proposition 2.4 it follows that the set \ $\ \{(x,y) : x \in [M, I), \ y \in
[M, I)^c, \ x \leq y\}\ $\  has positive product measure. But this is absurd
since the set is empty.
\end{mproof}

The interpolation property of the following theorem is quite crucial to
our proof of the main theorem
of section 4 and it will also be used in the proof of the
Kunneth formula. Note that if \ $\ x, y\ $\  belong
to \ $\ M ({\cal L})\ $\
and \ $\ x \leq y\ $\  then it is not necessarily true that
there exist
neighbourhoods \ $\ U_x\ $\  and \ $\ U_y\ $\
of \ $\ x,y\ $\  with \ $\ U_x \prec U_y\ $\ .
For this reason it seems that the proof of Theorem 2.6 must, inevitably,
be somewhat nontrivial.
\vspace{.3in}

\noindent {\bf Theorem 2.6}\it \ \  Let \ $\ x, y\ $\  be
distinct points of \ $\ M ({\cal L})\ $\  with
neighbourhoods \ $\ U_x, U_y\ $\
with \ $\ U_x \prec U_y.\ $\ Then there exist a point \ $\ z\ $\
in \ $\ M({\cal L})\ $\  and
neighbourhoods \ $\ V_x, V_y, V_z\ $\  of \ $\ x, y, z\ $\  such
that \ $\ V_x \prec V_z \prec V_y\ $\ . Furthermore, \ $z$\  can be chosen
in $\ U_x \ $, and $\ V_x $\ and $\  V_z$\ can be chosen to be subsets
of $\ U_x \ $.\rm

\begin{mproof}\  We may assume
that \ $\ U_x = [x, I) \cap C, \
U_y = [y, I) \cap C^\prime\ $\  with \ $\ C, C^\prime\ $\  as
in Proposition 2.4. If \ $\ U_x\ $\  is the
singleton \ $\ \{x\}\ $\  then we may simply
put \ $\ z =x,\  V_x = V_z = U_x, \ V_y = U_y\ $\ .

Suppose that \ $\ U_x\ $\  is not
a singleton and let \ $\ {\cal L}_1\ $\  be the
completely distributive projection
lattice \ $\ \tilde{P}(U_x){\cal L},\ $\  with
\ $\ 0_1\ $\  and \ $\ I_1\ $\  as the minimal and maximal projections. By the
hypothesis, and since \ $\cl$\ is completely distributive,
there exists a projection \ $\ F \in {\cl}_1\ $\
such that \ $\ F_- \neq I_1\ $\ .
The subset \ $\ [F_-, I_1)^c\ $\  of \ $\ M({\cal L}_1)\ $\  is nonvoid and
so, by Theorem 2.2,  there
is a projection \ $\ M_1\ $\  in \ $\ M({\cal L}_1)\ $\
with \ $\ F_- \leq M_1\ $.
Note that from the definition of \ $\ F_-\ $\  it follows
that \ $\ [F, I_1)^c \prec [F_-,I_1)\ $\
and so $\ [F,I_1)^c \prec [M_1,I_1)\ $
We now lift this comparability to \ $\ M({\cal L})\ $\  to get the desired new
neighbourhoods of the points \ $\ x\ $\  (the lift of \ $\ 0_1)\ $\
and \ $\ z\ $\  (the lift
of \ $\ M_1$).

This lifting is obtained with the map \ $\ \Phi : M ({\cal L}_1)
\rightarrow M({\cal L})\ $\  given by

\[
\Phi (M) = \sup \{L \in {\cal L} : L \tilde{P} (U_x) = M\}.
\]

\noindent Observe that  the set \ $\ \Phi ([F, I_1)^c)\ $\  is
equal to \ $\ [\Phi (F), I)^c \cap U_x.\ $\
To see this let \ $\ G \in [F, I_1)^c.\ $\  Then the projection
\ $\ K = G^\perp F \ $\  is a nonzero
projection in \ $\ {\cal L}_1.\ $\ If \ $\ L \tilde{P} (U_x) = G\ $\
then \ $\ L^\perp
\tilde{P}(U_x) \Phi (F) = L^\perp  \tilde{P}(U_x)F = G^\perp F = K,\ $\  and so
\ $\ \Phi(G)^\perp \tilde{P} (U_x) \Phi (F) = K,\ $\  and so \ $\ \Phi (G)
\not\geq \Phi (F).\ $\

Thus we now have that the set \ $\ V_x = \Phi ([F, I_1 )^c)\ $\  is an
open set in \ $M(\cl)$\ containing \ $\ x,\ $\  and is thus
a neighbourhood of \ $\ x.\ $\ Since
\ $\ \Phi\ $\  is order preserving, \ $\ V_x \prec V_z\ $\
where \ $\ V_z = \Phi ([F_-, I_1))\ $\
is a neighbourhood of \ $\ \Phi (M_1) = z\ $\  say. Since \ $\ V_z
\subseteq U_x\ $\  we
have \ $\ V_z \prec U_y\ $\  and so, with \ $\ V_y = U_y,\ $\
the proof is complete.
\end{mproof}

\section{Spectral Homology}

First, we define a point-neighbourhood homology for partially ordered
topological spaces.

Let $(X, \leq)$ be a separable topological space with
an antisymmetric partial order. An
{\em edge} of $(X, \leq)$, or {\em topological edge},
is an ordered pair $(x, y)$ of points of $X$ for
which there exist neighbourhoods $U_x, U_y$ such that $s \leq t$ for all $s \in
U_x$ and $t \in U_y$. If $F \subseteq X$ is a finite subset define the digraph
$G(F)$ to have vertex set $F$ and directed edges $(x, y)$ where $x, y$
belong to $F$ and
where $(x, y)$ is a topological edge of $(X, \leq)$. Associated with the
undirected graph of $G(F)$ is a simplicial
complex $\Delta(F)$ in which vertices
correspond to 0-simplices, edges to 1-simplices, and where complete
subgraphs on $t$ vertices correspond to $t-1$ simplices. Define $H_n
(\Delta(F))$ to be the usual integral simplicial homology of $\Delta(F)$ and
note that if $F \subseteq G$ are finite sets then the inclusion $\Delta(F)
\subseteq \Delta(G)$ gives a group
homomorphism $H_t (\Delta(F)) \rightarrow H_t (\Delta (G))$
for each $t$.
\vspace{.3in}

\noindent {\bf Definition 3.1} \ \  The {\it point-neighbourhood homology}
of $(X, \leq)$ are the
groups \ \ $H_t (X) = \displaystyle{\lim_\to}
H_t (\Delta (F))$, $t = 0, 1, \dots, $
where the direct
limit is taken over the net of finite subsets $\ F\ $of \ $X\ $.
\vspace{.3in}

In other words, $\ \ H_t(X)$\ \ is the integral simplicial
homology of the infinite simplicial complex, $\ \ \gD_{top}(X,\le)$\ \ say,
arising from the topological edges of $\ \ (X,\le).$\ \
We refer to $\ \ \gD_{top}(X,\le)$\ \ as the {\em topological complex}
of $\ \ (X,\le).$\ \
\vspace{.3in}

\noindent {\bf Definition 3.2} \ \ The {\it spectral homology} of
the reflexive operator algebra
${\cal A}$
associated with a completely distributive commutative projection lattice
${\cal L}$ is
defined to be the point-neighbourhood homology $H_*(M({\cal L}))$
of the partialy ordered set
$M({\cal L})$ of meet-irreducible projections
(excluding the identity operator).
The goups $H_t (M({\cal L}))$ are also written as
$H^{sp}_t ({\cal L})$ and $H^{sp}_t ({\cal A})$.

Write $\ \ H^{sp}_t ({\cal A}) = H^{sp}_t ({\cal L})=  H_t (M({\cal L}))$,\ \
for $\ \ t \ge 0$\ \ .\
\vspace{.3in}

We have chosen to define $H_*^{sp}(\ca)$ in terms of the topological complex
$\gD_{top}(M(\cl),\le)$ rather than the natural complex
$\gD (M(\cl),\le)$ of the partially ordered set $M(\cl)$.
The reason for this is that
this homology is adequate for the proofs of Theorem 4.1 and 4.2,
and it arises naturally in the arguments there involving comparability of
interval projections. However, in all of the examples in
this section the inclusion
$\gD_{top}(M(\cl),\le) \to \gD (M(\cl),\le)$, when it is proper,
nevertheless induces an isomorphism
of homology, and so it would be interesting to know if this
is a general phenomenon. (The proof of Theorem 3.9 provides a
little affirmative evidence.)
\vspace{.3in}

\noindent {\bf Proposition 3.3}\it $\ \
H_0^{sp} ({\cal A}) = \ZZ^d$ where $d$ is the linear space
dimension of the commutant of ${\cal A}$.\rm
\vspace{.3in}

\noindent {\bf Proposition 3.4} \it
\ \ If ${\cal A}$ is a nest
algebra then $H_t^{sp} ({\cal A}) = 0$ for all $t > 0$.\rm
\vspace{.3in}

\noindent {\bf Proposition 3.5}\em \ \ Spectral homology
is stable in the sense that
\[ H_t^{sp}
({\cal A}
\otimes {\cal L} ({\cal H})) = H^{sp}_t ({\cal A} \otimes M_k (\IC)) = H_t^{sp}
({\cal A})
\]
for all $t \geq 0$.\rm
\vspace{.3in}

\noindent {\bf Proposition 3.6}\it \ \ If ${\cal A} = A (G)$ is a
finite-dimensional digraph algebra,
associated with the digraph $G$, then, for $t \geq 0$,
$H^{sp}_t ({\cal A}) = H_t (\Delta
(G_r))$, the simplicial homology of the simplicial complex of the reduced graph
$G_r$ of $G$.\rm
\vspace{.3in}

\noindent {\bf Propositon 3.7}\it \ \ If ${\cal A}$ is an  infinite
tridiagonal algebra,
as in the introduction,
then $H_t^{sp} ({\cal A}) = 0$, for all \ $t > 0$.\rm
\vspace{.3in}

\noindent {\bf Proposition 3.8}\it \ \ Let $S{\cal A}$
denote the suspension of the CDCSL algebra
${\cal A}$
which is given by

\[
S{\cal A} = \left[ \begin{array}{ccc} {\cal B}_1 & 0 & *\\ 0 &
{\cal B}_2 & * \\ 0 & 0 &
{\cal A}
\end{array} \right]
\]

\noindent where ${\cal B}_1$ and ${\cal B}_2$ are type I factors.
Then $H_{t+1}^{sp}
(S{\cal A}) = H_t^{sp}({\cal A})$ for $t
\geq 1$.\rm
\vspace{.3in}

\noindent {\it Proofs:} (3.3) This is a consequence of Proposition 2.1 (iii).

(3.4) By Proposition 2.1 (i) $M({\cal L})$ is a totally ordered set when
${\cal L} =
$Lat$ {\cal A}$ is a nest. Also the topology on $M({\cal L})$ is the order
topology and so the proposition follows.

(3.5) Immediate from the fact that \ Lat(${\cal A}
\otimes {\cal L}({\cal H}))$ is
isomorphic to Lat
${\cal A}$.

(3.6) Let $\{e_{xy} : (x, y) \in E (G)\}$ be a (partial) system of matrix units
for ${\cal A}$, indexed by the edges of $G$. Define

\[
L_x = \sup \{L \in { Lat } {\cal A}: L e_{xx} = 0\}
\]

\noindent and verify that $x \rightarrow L_x$ effects an
isomorphism between $G_r = (V
(G_r), E(G_r))$ and $(M({\cal L}), \leq )$.
(The reduced graph of $\ \ G\ \ $
is, roughly speaking, the
antisymmetric graph obtained by collapsing each maximal complete subgraph
of equivalent vertices
to a single vertex.)
In fact the proposition is also
valid, with the same proof, for infinite digraph algebras - more usually
refered to as purely atomic CSL algebras.

(3.7) $M$(Lat$ {\cal A})$, as a graph, is an infinite (connected) chain whose
consecutive proper edges have alternating direction. There is one end in the
usual block staircase case, and no ends if the staircase is two-way infinite.
Plainly such chains have trivial simplicial homology groups for \ $t > 0$.

(3.8) Let $L \in M ($Lat$ {\cal A})$.
Then $\tilde{L} = I \oplus I \oplus L$ belongs to
$M($Lat$ (S{\cal A}))$.
Also the projections $P = I \oplus 0 \oplus I$ and $Q = 0 \oplus I \oplus I$
belong to $M($Lat$ ({\cal A}))$, and $M($Lat$ (S{\cal A})) = \{P, Q\}
\cup \{\tilde{L} : L \in M
(Lat {\cal A})\}$. This partially ordered set is isomorphic to the
infinite digraph
arising from the two point
suspension $S(M({\cal L}  ))$ of the infinite digraph
$M({\cal L})$.
Since $H_{t+1} (S(G)) = H_t (G)$ for finite
digraphs (for $t \geq 1$), by elementary simplical
homology, the conclusion follows readily. \hfill $\Box$
\vspace{.3in}

\noindent {\bf Theorem 3.9}\it  \ \ The spectral homology of the
spatial tensor product ${\cal A} \otimes {\cal A}^\prime$ of two
CDCSL algebras is computable by the Kunneth formula:

\[
H_n^{sp}({\cal A} \otimes {\cal A}^\prime) = (\sum_{p+q=n} \oplus \
(H_p^{sp}({\cal A})
\otimes H_q^{sp}({\cal A}^\prime)) \ \ )\oplus \ \ ( \sum_{p+q=n-1}
 \oplus \ {\rm Tor}(H_p^{sp}({\cal A}),H_q^{sp}({\cal A}^\prime))).
\]
\rm
\begin{mproof}
First, we identify the product structure of
\ $\ M(\cl_1 \times \cl_2)\ $,\  and we begin by showing that
if $\ \ E \in \cl_1 \ot \cl_2 = \cl$\ \  then $\ \ E_+$\ \  has the
form $\ \ F_+ \ot G_+$\ \  for
some $\ \ F$\ \  in $\ \ \cl_1$\ \  and $\ \ G$\ \  in $\ \ \cl_2$.\ \
A lattice theoretic approach to this product structure can be found in
Fraser \cite{fra}.

Choose $\ \ F_\gb \in \cl_1$\ \  and $\ \ G_\gb \in \cl_2$,\ \
for $\ \ \gb$\ \  in some index
set, so that $\ \ E^\bot = $\ \ sup$\{F_\gb^\bot \ot G_\gb^\bot\}$.\ \
That this is possible follows
from the fact that the complementary lattice of $\ \ \cl$,\ \  the lattice
of projections $\ \ L^\bot$\ \  with $\ \ L \in \cl$,\ \  is the spatial
tensor product of the complementary lattices
of $\ \ \cl_1$\ \  and $\ \ \cl_2$.\ \
Set $\ \ F= $\ inf$\{F_\gb\}, G= $\ inf$\{G_\gb\}.$\ \
If $\ \ H \in \cl_1$\ \  and $\ \ H \not\le F$\ \
 then $\ \ H \not\le F_\gb$\ \  for some index $\ \ \gb$\ \
and so $\ \ HF_\gb^\bot \ne 0$.\ \  Thus
$\ \ (H \ot I)(F_\gb^\bot \ot G_\gb^\bot) \ne 0$\ \ ,
and so $\ \ (H \ot I)E^\bot \ne 0$.\ \  That is, $\ \ H \ot I \not\le E$.\ \
Taking the infimum over such $\ \ H$\ \  obtain $\ \ F_+ \ot I \ge E_+$.\ \
Similarly $\ \ I \ot G_+ \ge E_+$,\ \  and so $\ \ F_+ \ot G_+ \ge E_+$.\ \

On the other hand, let $\ \ K \not\le E $\ \  and express $\ \ K^\bot$\ \  as
\ \ sup$\{L_\ga^\bot \ot M_\ga^\bot\}$\ \
with $\ \ \ga$\ \  ranging
over some index set. Since $\ \ K^\bot \not\ge E^\bot$\ \  it follows that
$\ \ F_\gb^\bot \ot G_\gb^\bot \not\le K^\bot$\ \
for some index $\ \ \gb$.\ \  Thus,
for each index $\ \ \ga$\ \  either
$\ \ F_\gb^\bot \not\le L_\ga^\bot$\ \  or
$\ \ G_\gb^\bot \not\le M_\ga^\bot$,\ \  that is, either
$\ \ L_\ga \not\le F_\gb$\ \
or $\ \ M_\ga \not\le G_\gb$\ \
for each $\ \ \ga$.\ \
Thus,
\[
K^\bot  = (\bigvee_{L_\ga \not\le F_\gb} L_\ga^\bot \ot M_\gb^\bot) \vee
(\bigvee_{M_\ga \not\le G_\gb} L_\ga^\bot \ot M_\gb^\bot)
\]\[
\le (\bigvee_{L_\ga \not\le F_\gb} (L_\ga \ot I)^\bot) \vee
(\bigvee_{M_\ga \not\le G_\gb} (I \ot M_\ga)^\bot)
\]\[
= (\bigwedge_{L_\ga \not\le F_\gb} L_\ga \ot I)^\bot) \vee
(\bigwedge_{M_\ga \not\le G_\gb} I \ot M_\ga)^\bot
\]\[
\le ((F_\gb)_+ \ot I)^\bot \vee (I \ot (G_\gb)_+)^\bot
\]

\noindent and so $\ \ K \ge (F_\gb)_+ \ot (G_\gb)_+ \ge F_+ \ot G_+.$\ \

Now let $\ \ L \in M(\cl_1 \ot \cl_2)$\ \  and set
\[
M_1 = \mbox{sup}\{L_1 \in \cl_1: L_1 \ot I \le L\}
\]\[
M_2 = \mbox{sup}\{L_2 \in \cl_2: I \ot L_2  \le L\}
\]
so that, with tolerable abuse of
notation,  $\ \ M_1 \vee M_2 \le L$.\ \  Assume, by way of contradiction,
that $\ \ M_1 \vee M_2 \ne L.$\ \
Since $\ \ \cl_1 \ot \cl_2$\ \  is completely distributive
there exists a projection $\ \ E$\ \  in
$\ \ \cl_1 \ot \cl_2$\ \
such that $\ \ E_+(M_1 \vee M_2)^\bot \ne 0$\ \  and $\ \ E_+L^\bot = 0.$\ \
{}From the above, $\ \ E_+ = F_+ \ot G_+$\ \
for some $\ \ F \in \cl_1$\ \  and $\ \ G$\ \  in $\ \ \cl_2$,\ \
and so $\ \ 0 \ne (F_+ \ot G_+)(M_1 \vee M_2)^\bot =
 (F_+ \ot G_+)(M_1^\bot \ot M_2^\bot) = F_+M_1^\bot \ot G_+M_2^\bot.$\ \
Thus $\ \ 0 \ne F_+M_1^\bot$\ \  and $\ \ 0 \ne G_+M_2^\bot$.\ \  But now
$\ \ L = (M_1 + F_+ M_1^\bot)(M_2 + G_+M_2^\bot), $\ \
and so $\ \ L$\ \  is not meet-irreducible, contrary to our assumption.

Finally, it is straightforward to show that if $\ \ L_i \in M(\cl_i)$,\ \
for $\ \ i = 1,2,$\ \
then $\ \ L_1 \vee L_2$\ \  is an element of $\ \ M(\cl_1 \ot \cl_2)$.\ \

We may now identify  $\ M(\cl_1 \otimes \cl_2)\ $ and
its ordering with the set
$\ M(\cl_1) \times M(\cl_2)\ $ with the product order.
In particular the complex
$\gD(M(L),\le)$
coincides with the usual product complex
$\ \gD(M(L_1),\le) \times \gD(M(L_2),\le)$.\ \
On the
other hand, because of the definition of the product
topology, one sees that $\ \gD_{top}(M(L))$\ is naturally
isomorphic to what we call the edge product complex
$\ \ \gD_{top}(M(L_1)) \times_e \gD_{top}(M(L_2)).$\ \
This is defined to be the complex which is determined by the
digraph with the product vertex set $\ \ M(\cl_1) \times M(\cl_2)$\ \
and edges $\ ((x_1,y_1),(x_2,y_2))$\
arising from the topological
edges $\ (x_1,x_2),\  (y_1,y_2)$  of $\ M(\cl_1)$\
and $\ M(\cl_2)$\  respectively. Nevertheless
we shall now show that the inclusion
of \ \ $\gD_{top}(M(L_1)) \times_e \gD_{top}(M(L_2))$\ \ in
$\ \ \gD_{top}(M(L_1)) \times \gD_{top}(M(L_2))$\ \
induces an isomorphism of  homology. This will imply that
\[
 H^{sp}_n(\ca \otimes \ca^\prime) =
H^{sp}_n(M(\cl) =
\]
\[
H_n(\gD_{top}(M(L_1)) \times_e \gD_{top}(M(L_2))) =
H_n(\gD_{top}(M(L_1)) \times \gD_{top}(M(L_2))),
\]
and the desired Kunnneth formula then follows from the
corresponding formula  in simplicial homology.

To this end let \ $\ F_i \subseteq \cl_i\ $\  be finite sets
determining the subcomplexes
$\ \gD(F_i)$\ of $\ \gD_{top}(M(L_i).$\

Each
1-simplex of
$\ \gD(F_1) \times \gD(F_2)\ $
which is not a 1-simplex of
$\ \gD_{top}(M(L),\le)$
comes from an edge of the form
\ $\ ((x,y_1),(x,y_2))\ $\  or \ $\ ((x_1,y),(x_2,y))\ $,\  where \ $\ x\ $\
(resp. \ $\ y\ $\ )
is a point of \ $\ F_1\ $\
(resp. \ $\ F_2\ $\ ) such that \ $\ \{x\}\ $\  (resp. \ $\ \{y\}\ $\ ) is
not a neighbourhood of \ $\ x\ $\
(resp. \ $\ y\ $\ ). Refer to such edges, and their simplices,
as the {\em extremal} edges and simplices of
\ $\ \gD(F_1) \times \gD(F_2)\ $.\
In the case of the extremal edge
\ $\ ((x_1,y),(x_2,y))\ $,\ use Theorem 2.6 to
choose a point \ $\ y^\prime\ $\  in a neighbourhood of \ $\ y\ $\
and choose \ $\ x_2^\prime\ $\  in a neighbourhood of \ $\ x_2\ $\  so
that \ $\ (y,y^\prime)\ $\  and
\ $\ (x_2,x_2^\prime)\ $\  are topological edges. Then
\ $\ ((x_1,y),(x_2^\prime ,y^\prime))\ $\  and \ $\ ((x_2,y),(x_2^\prime,
y^\prime))\ $\  are topological edges. Assume, moreover, that the
point \ $\ (x_2,y)\ $\  in
\ $\ F_1 \times F_2\ $\  is maximal in the partial ordering, so that each
1-simplex of $\ \gD(F_1) \times \gD(F_2)\ $
,\  with \ $\ (x,y_2)\ $\  as an endpoint,
arises from an edge of the form above.
By choosing \ $\ (x_2^\prime ,y^\prime)\ $\  in a
sufficiently small neighbourhood we can
ensure that whenever
 \ $\ ((u,v),(x_2,y))\ $\  is an edge corresponding to a 1-simplex of
\ $\ \gD(F_1) \times \gD(F_2)\ $\ (even an extremal edge) then
\ $\ ((u,v),(x_2^\prime,y^\prime))\ $\  is a topological edge for
\ $\ M(\cl)\ $.\  The product complex
\ $\ \gD(F_1) \times \gD(F_2)\ $\  is now chain homotopic
to a subcomplex \ $\ \gD\ $\
of \ $\ \gD(F_1 \cup \{x_2^\prime\}) \times \gD(F_2 \cup \{y^\prime\})\ $\
through the chain
homotopy which moves \ $\ (x_2,y)\ $\  to \ $\ (x_2^\prime ,y^\prime)\ $.\
The point of this
homotopy is that the complex
\ $\ \gD\ $\  has fewer extremal edges. Note that maximal
points of the type above always
exist. Thus, repeating such homotopies a finite number of times,
we find that
there exist finite sets \ $\ G_1\ $\  and \ $\ G_2\ $,\
containing \ $\ F_1\ $\  and \ $\ F_2\ $,\  respectively,
so that  the product complex \ $\ \gD(F_1) \times \gD(F_2)\ $\  is
homotopic in \ $\ \gD(G_1) \times \gD(G_2)\ $\  to a complex with no extremal
edges. That is,
\ $\ \gD(F_1) \times \gD(F_2)\ $\  is homotopic to a subcomplex
of $\ \gD_{top}(M(\cl),\le)$,\
and the desired homology isomorphism follows.
\end{mproof}
\vspace{.3in}

\noindent {\bf Examples 3.10 (a)}
\ \  If \ $\ \ca\ $\  is a tridiagonal algebra, as in the introduction,
then $H_1^{sp}(T_{\IR} \otimes \ca)) = 0$.

\noindent {\bf (b)}
\ \ Let \ $\ \ct_{[0,1]}\ $\  be the Volterra nest
algebra, that is, the nest algebra on \ $\ L^2[0,1]\ $\  associated with the
nest of intervals \ $\ [0,t], \ 0 \le t \le 1.\ $\ Then \ $\
H^{sp}_1(\ct_{[0,1]} \ot A(D_4)) = \ZZ.\ $\
This a consequence of Propositions 3.6 and Theorem 3.9. However it is
also posssible to obtain this by direct methods, as follows.

First, show directly that \ $\ (M(\cl),\le)\ $\  is isomorphic to the product
set \ $\ X = [0,1] \times \{1,2,3,4\}\ $\  with the natural product order :
\ $\ (x,n) \le (y,m)\ $\  if and only if \ $\ x \le y\ $\  and \ $\ (n,m) \ $\
is an edge
of the 4-cycle digraph \ $\ D_4.\ $\  Because of the nature
of the product topology
the topological edges are the pairs \ $\ ((x,n),(y,m))\ $\  for which \ $\ x <
y\ $\  and
\ $\ (n,m)\ $\  is an edge of \ $\ D_4.\ $\   Next, let $F = F_0 \times
\{1,2,3,4\}
\ $\  where \ $\ F_o\ $\  is a finite subset
of \ $\ [0,1].\ $\   It is an elemenary exercise in simplicial homology to show
that
the subcomplex \ $\ \gD(F)\ $\  of the topological complex \ $\ \gD(X)\ $\
has first simplicial homology group equal to \ $\ \ZZ.\ $\
Furthermore the inclusion \ $\ \gD(F) \to \gD(X)\ $\
induces the identity map \ $\ \ZZ \to \ZZ.\ $\  Since every
finite subcomplex of \ $\ \gD(X)\ $\ is contained in one of these special
subcomplexes
it follows that \ $\ H_1^{sp}(\cl) =\ZZ.\ $\

\noindent  {\bf (c)} \ \
 Let $G$ be an antisymmetric digraph  for which $| \Delta (G)|$, the
geometric realisation of the complex of $G$, is
the Klein bottle. Then $H_1 (\ct_{[0,1]} \otimes A(G)) = \ZZ_2$, and $H_1
(\ct_{[0,1]}
\otimes A(D_4) \otimes A(G)) = \ZZ_2 \oplus \ZZ$.

\noindent {\bf (d)}.
Let $\ X = [-1,1]^2$\   be the closed square in $\ \IR^2$\
carrying the Borel measure $\ \mu$\  which
is the sum of Lebesgue area measure and arc length measure for the two
diagonals. There is a partial ordering of $\ X$\  which is suggested by the
following diagram:

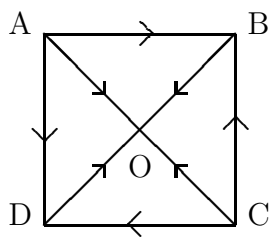
\begin{center}
\setlength{\unitlength}{0.0125in}%
\begin{picture}(100,94)(65,715)
\thicklines
\put( 80,720){\line( 1, 0){ 80}}
\multiput(120,795)(0.41667,0.41667){13}{\makebox(0.4444,0.6667){\sevrm .}}
\multiput(125,800)(-0.41667,0.41667){13}{\makebox(0.4444,0.6667){\sevrm .}}
\multiput( 75,760)(0.41667,-0.41667){13}{\makebox(0.4444,0.6667){\sevrm .}}
\multiput( 80,755)(0.41667,0.41667){13}{\makebox(0.4444,0.6667){\sevrm .}}
\multiput(120,725)(-0.41667,-0.41667){13}{\makebox(0.4444,0.6667){\sevrm .}}
\multiput(115,720)(0.41667,-0.41667){13}{\makebox(0.4444,0.6667){\sevrm .}}
\multiput(165,760)(-0.41667,0.41667){13}{\makebox(0.4444,0.6667){\sevrm .}}
\multiput(160,765)(-0.41667,-0.41667){13}{\makebox(0.4444,0.6667){\sevrm .}}
\put(105,780){\line( 0,-1){  5}}
\put(105,775){\line(-1, 0){  5}}
\put( 80,800){\line( 1, 0){ 80}}
\put(160,800){\line( 0,-1){ 80}}
\put(160,720){\line(-1, 1){ 80}}
\put( 80,800){\line( 0,-1){ 80}}
\put( 80,720){\line( 1, 1){ 80}}
\put(100,745){\line( 1, 0){  5}}
\put(105,745){\line( 0,-1){  5}}
\put(115,740){\makebox(0,0)[lb]{\raisebox{0pt}[0pt][0pt]{\twlrm O}}}
\put(135,780){\line( 0,-1){  5}}
\put(135,775){\line( 1, 0){  5}}
\put(140,745){\line(-1, 0){  5}}
\put(135,745){\line( 0,-1){  5}}
\put( 65,800){\makebox(0,0)[lb]{\raisebox{0pt}[0pt][0pt]{\twlrm A}}}
\put(165,800){\makebox(0,0)[lb]{\raisebox{0pt}[0pt][0pt]{\twlrm B}}}
\put( 65,720){\makebox(0,0)[lb]{\raisebox{0pt}[0pt][0pt]{\twlrm D}}}
\put(165,720){\makebox(0,0)[lb]{\raisebox{0pt}[0pt][0pt]{\twlrm C}}}
\end{picture}
\end{center}

\noindent  That is, if  points $\ P, Q $\  belong to a triangle, say the
triangle of $\ ABO,$\
then $\ Q \le P$\  if and only if
the vector $\ PQ $\ is equal to the vector $\  sAB + tBO$\  for some $\ s, t
\ge 0.$\
The triangles $\ CBO, CDO, ADO$\  carry
similar such orders, which are compatible on overlapping edges, and the
partial ordering on $\ X$\  is the union of these four partial orderings. The
lattice of decreasing Borel sets for $\ (X, \le, \mu)$\
gives a completely distributive
commutative lattice $\ \cl$\  of projections on $\ L^2(X,\mu),$\
and an associated reflexive
operator algebra $\ \ca.$\
  It may seem curious, at first glance, that $\ H^{sp}_{1}(\ca) = \ZZ,$\
since there is no apparent
hole to contribute 1-cycles that are not 1-boundaries.
But the origin has measure
zero, so it may be deleted, and the assertion now becomes plausible ! In
fact examination reveals that $\ (M(\cl), \le)$\
is naturally isomorphic to the set
\[
\{(-1,1)^2 \backslash \{(0,0)\}\} \cup \{(-1,1),(1,1),(-1,-1),(1,-1)\}
\]
together with the relative ordering from $\ [-1,1]^2,$\
already described, and the first point-neighbourhood homology group of this
partially ordered set is $\ \ZZ.$\
\vspace{.3in}

Gilfeather and Smith \cite{gil-smi-2} have obtained a Kunneth style formula
for the Hochschild cohomology of the join of two operator
algebras, $\ \ \ca \# \cb\ \ $. This is
shown to be valid if one of the algebras acts on a finite-
dimensional Hilbert space, and is shown to be false in general.
For the spectral homology of CDCSL algebras the situation is much
more straightforward. Also, we see that this context is considerably
simpler than
that of Theorem 3.9.
\vspace{.3in}

\noindent {\bf Proposition 3.11}\it \ \ Let $\ \ \ca\ \ $ and $\ \ \cb\ \ $
be CDCSL algebras and let $\ \ \ca \# \cb\ \ $ be their join :
\[
\ca \# \cb = \left[ \begin{array}{cc}
\cb & 0 \\
\star &{\cal A}
\end{array} \right].
\]
Then
\[
H_n(\ca \# \cb) = (\sum_{p+q=n-1} \oplus \ (H_p({\cal A})
\otimes H_q(\cb))).
\]
\rm

\begin{mproof}.
It is elementary to verify that $\ \ M($Lat$(\ca \# \cb))$
is the set of projections of the form
$\ \ M \oplus I\ \ $ or \ \ $0 \oplus L\ \ $ where $\ \ L\ \ $ belongs to
\ \ $M($Lat$(\ca))$\ \  and $\ \ M\ \ $ belongs to
\ \ $M($Lat$(\cb))$.\ \  It follows that the
 topological complex for the join algebra is precisely the simplicial complex
join of the topological complexes for $\ \ \ca\ \ $ and $\ \ \cb.\ \ $ Thus the
desired formula follows from the corresponding formula in simplicial homology.
\end{mproof}

We have seen that various constructions in simplicial homology,
such as joins, suspensions and products, are also available
for partially ordered measure spaces, and
for CSL algebras.
Similarly, the following somewhat unusual CSL algebra, a fibre sum
of nest algebras, can be defined at the algebraic level, as a pull back,
as below,
or in terms of a fibre sum of the constituent partially ordered
measure spaces. Such constructions can be used to create algebras, such as
\ $(d)$\ above, with interesting homology.
\vspace{.3in}

\noindent {\bf Example 3.12}\ \
Let \ $\ \ct_{\mu}\ $\  be the nest algebra on \ $\ L^2([0,1],\mu)\ $\
where \ $\ \mu\ $\  is equal to
Lebesgue measure plus unit masses at \ $\ 0\ $\  and at \ $\ 1.\ $\
Let \ $\ \ct_{\mu} \op_{\ \IC} \ct_{\mu}\ $\
be the fibre sum algebra associated with the summand maps
\ $\ \ct_{\mu} \to \IC\ $\  given by  the compression
maps for the atomic interval
projections over 1. Let
\[
\ca = (\ct_{\mu} \op_{\ \IC} \ct_{\mu}) \oplus_{\ \IC^2} (\ct_{\mu} \op_{\ \IC}
\ct_{\mu})
\]
be the fibre sum algebra arising from summand maps
\[
\ct_{\mu} \op_{\ \IC} \ct_{\mu} \to \IC^2
\]
corresponding to the compression maps for the  atomic
intervals associated with \ $\ 0\ $\  in each copy of \ $\ [0,1].\ $\
Then \ $\ \ca\ $\  is a completely distributive CSL algebra for which \ $\
M(\cl),\ $\
as a set,
consists of  four unit intervals joined at their endpoints to form a square.
The
partial ordering is that inherited from the four intervals and
on consideration of the (infinite) topological complex of $M(\cl)$\
it follows
that
\ $\ H_1^{sp}(\ca) = \ZZ.\ $\

\section{The Main Result}

Let \ \ ${\cal A}\ \ $ be a weak star closed operator algebra. Let \ \ $Z^1
({\cal A})\ \ $ denote the space
of derivations that are weak star continuous and let \ \ $B^1_{ess}({\cal A})\
\ $ denote
the subspace of such derivations which are point weak star limits of inner
derivations.  Thus \ \ $\delta \in B^1_{ess}({\cal A})\ \ $ if and only
if for each
weak star topology neighbourhood $\ \ U \subseteq \ca$\ \ of zero
and for each finite set \ \ $a_1,\dots,a_n\ \ $ in \ \ $\ca\ \ $ there exists
an
inner derivation \ \ $\delta_0\ \ $ such that
\ \ $ \delta(a_i) - \delta_0(a_i) \in U \ \ $ for \ \ $1 \le i \le n.\ \ $
We define the {\it essential Hochschild cohomology} group  Hoch$^1_{ess} ({\cal
A}) \ \ $ to be the
space \ \ $Z^1({\cal A}) / B^1_{ess} ({\cal A}).\ \ $

A Schur automorphism of \ \ ${\cal A}\ \ $ is, by definition, an automorphism \
\ $\alpha\ \ $ for
which there exists a masa \ \ ${\cal C}\ \ $  in \ \ ${\cal A}\ \ $ which is
fixed elementwise by \ \ $\alpha.\ \ $
\vspace{.3in}

\noindent {\bf  Theorem 4.1}\it \ \ Let \ \ ${\cal A}\ \ $ be a CDCSL algebra
for which \ \ $H_1^{sp} ({\cal A}) = 0.\ \ $ Then
Hoch$_{ess}^1 ({\cal A}) = 0.\ \ $
\vspace{.3in}\rm

This theorem will be a consequence of the following closely related  result.
\vspace{.3in}

\noindent {\bf Theorem 4.2}\it \ \ Let \ \ ${\cal A}\ \ $ be a CDCSL algebra
with \ \ $H^{sp}_1 ({\cal A}) = 0\ \ $ and let
\ \ $\alpha\ \ $ be a Schur automorphism of \ \ ${\cal A}\ \ $ with respect to
a masa \ \ ${\cal C}.\ \ $ Then
\ \ $\alpha\ \ $ is a point weak star limit of inner automorphisms of the form
\ \ $A
\rightarrow CAC^{-1}\ \ $ with \ \ $C \in \cc.\ \ $\rm
\vspace{.3in}

Fix \ \ ${\cal A}\ \ $ and \ \ ${\cal C}\ \ $ as in Theorem 4.2 and choose a
spectral representation for
\ \ ${\cal C}\ \ $ as follows. Identify the underlying
Hilbert space \ \ $H\ \ $ with \ \ $L^2 (m)\ \ $ for some finite measure space
\ \ $(Y,
\sigma, m)\ \ $ in such a way that \ \ ${\cal C}\ \ $ is identified with \ \
$L^\infty (m)\ \ $ by means of
multiplication operators. The lattice \ \ ${\cal L} = \ \ $Lat${\cal A} \ \ $
is a subset of \ \ ${\cal C},\ \ $ and we
shall frequently consider partitions of the identity by the atoms \ \ $Q_i\ \ $
of a
finite sublattice of \ \ ${\cal L}.\ \ $ Such partitions correspond to
particular finite
measurable partitions of \ \ $Y.\ \ $
\vspace{.3in}

\noindent {\bf Lemma 4.3}\it \ \  If \ \ $Q_1, \ Q_2, \ Q_3\ \ $ are interval
projections of \ \ ${\cal L}\ \ $ with \ \ $Q_1 \prec
Q_2 \prec Q_3\ \ $ then there exist invertible operators \ \ $C_1, \
C_3\ \ $ in \ \ ${\cal C}\ \ $ such that

\[
\alpha (Q_1 A Q_3) = C_1 Q_1 A Q_3 C_3
\]

\noindent for all \ \ $A\ \ $ in \ \ $L(H).\ \ $\rm

\begin{mproof}\  Consider the Hilbert space decomposition \ \ $H_1 \oplus \dots
\oplus H_4\ \ $
where \ \ $H_i = Q_i H, $\ for $\ 1 \leq i \leq 3,\ \ $ and \ \ $H_4 = Q_4 H\ \
$ with \ \ $Q_4 = I - Q_1 -
Q_2 - Q_3.\ \ $ Then

\[
\alpha: \left[ \begin{array}{cccc} 0 & X_1 & X_3 & 0 \\0 & 0 & X_2 & 0 \\ 0 & 0
& 0 & 0 \\ 0 & 0 & 0 & 0 \end{array} \right] \rightarrow \left[
\begin{array}{cccc} 0 & \alpha_1 (X_1) & \alpha_3 (X_3) & 0 \\ 0 & 0 & \alpha_2
(X_2) & 0 \\ 0 & 0 & 0 & 0\\ 0 & 0 & 0 & 0 \end{array} \right]
\]

\noindent where \ \ $\alpha_1, \alpha_2, \alpha_3\ \ $ are bimodule maps with
respect to the pairs
of masas \ \ $(Q_1 \cc, Q_2 \cc), \ (Q_2 \cc, Q_3 \cc)\ \ $ and \ \ $(Q_1 \cc,
Q_3 \cc)\ \ $
respectively. By a result of Haagerup  \cite{haa} each \ \ $\alpha_i\ \ $ is
completely bounded
and has the form

\[
\alpha_i: X \rightarrow \sum_k \phi_{k,i} X_i \psi_{k, i},
\]

\noindent where \ \ $\phi_{k,i} \ \ $ and \ \ $\psi_{k,i} \ \ $ belong to the
appropriate restriction of \ \ ${\cal C}\ \ $
and
\[
(\sum_k ||\phi_{k,i} ||^2 )^{\frac{1}{2}} (\sum_k || \psi_{k, i}
||^2 )^{\frac{1}{2}}\leq ||\alpha_i ||_{cb}.
\]

\noindent View the operators \ \ $\phi_{k,i}, \ \psi_{k, i}\ \ $ as elements in
\ \ $L^\infty (m).\ \ $
{}From the inequality we see that the function \ \ $\Phi_i (x, y) = \sum_k
\phi_{k, i}
(x) \psi_{k, i} (y)\ \ $ defines an element in \ \ $L^\infty (m \times m).\ \ $
In
particular the restriction of $\ \ \ga_i$\ \ to the subspace of
Hilbert-Schmidt operators
coincides with the
Schur multiplier of the representing kernel functions induced by $\ \ \Phi.$\ \

Since \ \ $\alpha\ \ $ is an automorphism we have \ \ $\alpha_1 (X_1) \alpha_2
(X_2) =
\alpha_3 (X_1 X_2).\ \ $ From this it follows, by considering the cases when \
\ $X_1,
X_2\ \ $ have rank one, for example, that for almost every pair \ \ $(x, z)\ \
$

\[
\Phi_1 (x, y) \Phi_2 (y, z) = \Phi_3 (x, z)
\]

\noindent for almost every \ \ $y.\ \ $ Thus we obtain a factorisation \ \
$\Phi_3 (x, y) = \theta (x)
\eta (y)\ \ $ with \ \ $\theta, \eta\ \ $ in \ \ $L^\infty (m)\ \ $ and hence
the representation
\ \ $\alpha_3 (X_3) = C_1 X_3 C_3,\ \ $ at least for Hilbert Schmidt operators.
But
\ \ $\alpha\ \ $ is weak star continuous, by Theorem 2.2 of \cite{dav-scp},
and so
this equality holds generally. Since
\ \ $\alpha_3\ \ $ is bounded below, and \ \ $C_1, C_3\ \ $ are bounded, it
follows that \ \ $C_1,
C_3\ \ $ are invertible.
\end{mproof}

Suppose now that \ \ $Q_1 \cl (H) Q_i \subseteq \ca,\ \ $ for \ \ $i = 2, 3,\ \
$ and that there exist \ \ $C_2
\in Q_1 \cc, D_2 \in Q_2 \cc, C_3 \in Q_1 \cc, D_3 \in Q_3 \cc,\ \ $ such that
these
operators are invertible in their respective spaces, and

\[
\alpha (Q_1 X Q_2) = C_2 Q_1 X Q_2 D_2,
\]
\[
\alpha (Q_1 Y Q_3) = C_3 Q_1 Y Q_3 D_3,
\]
for all \ \ $X, Y\ \ $ in \ \ $\cl(H).\ \ $ Then we claim that \ \ $C_2 =
\lambda C_3\ \ $ for some
nonzero scalar \ \ $\lambda.\ \ $ To see this we make use of the following
well-known
facts.
\vspace{.3in}

\noindent {\bf Lemma 4.4}\it \ \ Let \ \ ${\cal B}\ \ $ be a CDCSL
algebra with invariant projection lattice
\ \ ${\cal L}.\ \ $

{}~(i) If \ \ $ L \in \cl $\  and if $\ \ R = x  \otimes y^*\ \ $ is a nonzero
rank one operator with range vector
\ \ $x \in L\cal{H}\ \ $ and domain vector \ \ $y \in (I-L_-) \cal{H}\ \ $ then
\ \ $R \in {\cal B}.\ \ $

(ii) \ \ $\sup \{ {L \in  L} : I - L_- \neq 0\} = I.\ \ $
\vspace{.3in}\rm

To establish the claim note that the algebra  \ \ ${\cal B}_1 =
Q_1 {\cal A} Q_1\ \ $
is a CDCSL algebra on \ \ $Q_1 \cal{H}\ \ $
and let \ \ $R = x \otimes y^*\ \ $ be a nonzero rank one operator in \ \
${\cal B}_1.\ \ $ If \ \ $X
= Q_1 X Q_2\ \ $ and \ \ $Y = Q_1 YQ_3\ \ $ then, since \ \ $\alpha\ \ $ is an
automorphism,
the
following equations hold :

\[
C_2 R X D_2 = \alpha (R) C_2 X D_2,
\]
\[
C_3 RYD_3 = \alpha(R) C_3 YD_3
\]

\noindent From the first equation it follows that for all choices of \ \ $X\ \
$ the range of
\ \ $\alpha(R) C_2 X D_2\ \ $ is \ \ $\{\lambda C_2 x: \lambda \in \IC\}.\ \ $
Since \ \ $C_2\ \ $ and
\ \ $D_2\ \ $ are invertible on their respective spaces \ \ $\alpha (R)\ \ $ is
necessarily of
rank one. We remark at this point that an automorphism of a CDCSL algebra need
not preserve rank - see \cite{gil-moo}.
Considering both equations, obtain \ \ $C_2
x = \lambda_x C_3 x\ \ $ for some nonzero scalar \ \ $\lambda_x.\ \ $ The
obvious
manipulations show that if \ \ $x_1\ \ $ and \ \ $x_2\ \ $ are linearly
independent then
\ \ $\lambda_{x_1} = \lambda_{x_2},\ \ $ and if \ \ $y = \mu x\ \ $
then \ \ $\lambda_{y} =
\lambda_x.\ \ $ Thus, using the second part of the lemma, obtain \ \ $C_2 =
\lambda C_3\ \ $ for some scalar \ \ $\lambda,\ \ $ and the claim is proven.

Suppose now that \ \ $Q_1 \prec Q_2 \prec Q_3\ \ $ and it is known that for any
\ \ $X\ \ $ of the form
\ \ $X = Q_1 X Q_2\ \ $ and for any \ \ $Y\ \ $ of the form \ \ $Q_2 Y Q_3\ \ $
that \ \ $\alpha (X) =
C_1 X D_1,~ \alpha(Y) = C_2 Y D_2\ \ $ with \ \ $C_1 \in Q_1 \cc,~ D_1, C_2 \in
Q_2 \cc\ \ $
and \ \ $D_2 \in Q_3 \cc,\ \ $ where each operator is invertible in its
respective
space. Suppose further, that it is known that \ \ $\alpha (Z) = C_3 Z D_3\ \ $
when \ \ $Z =
Q_1 Z Q_3,\ \ $ with \ \ $C_3 \in Q_1 \cc,~ D_3 \in Q_3 \cc\ \ $ invertible.
Then
since \ \ $\alpha\ \ $ is an automorphism,

\[
C_3 XYD_3 = C_1 XD_1 C_2 YD_2.
\]

\noindent From this equality it follows readily that \ \ $C_2 = \lambda
D^{-1}_1\ \ $ for some nonzero
scalar \ \ $\lambda.\ \ $

Motivated by the relationships above between the various local implenting
operators we now formulate a general lemma which is exactly suited to our
needs.
\vspace{.3in}

\noindent {\bf Lemma 4.5}\it \ \ Let \ \ $Y = \{ x_1, \dots, x_m\}\ \ $ be a
finite set with an
antisymmetric partial ordering. Let \ \ $G_1, \dots, G_m\ \ $ be abelian
groups with \ \ $\IC_* \subseteq  G_i\ \ $ for all \ \ $i\ \ $ and suppose that
for each pair \ \ $x_i
\prec
x_j\ \ $ there is associated a pair \ \ $(c_{ij}, d_{ij})\ \ $ with \ \ $c_{ij}
\in G_i\ \ $ and
\ \ $d_{ij} \in G_j\ \ $ such that the following properties hold:

{}~~(i) If \ \ $x_i \prec x_j\ \ $ and \ \ $x_i \prec x_k\ \ $ then \ \ $c_{ij}
= \lambda c_{ik}\ \ $ for some
\ \ $\lambda \in \IC_*.\ \ $

{}~(ii) If \ \ $x_i \prec x_j\ \ $ and \ \ $x_k \prec x_j\ \ $ then \ \ $d_{ij}
= \lambda d_{kj}\ \ $ for some
\ \ $\lambda \in \IC_*.\ \ $

{}~(iii)  If \ \ $x_i \prec x_j \prec x_k\ \ $ then \ \ $d_{ij} = \lambda
c_{jk}\ \ $ for some scalar
\ \ $\lambda \in \IC_*.\ \ $

{}~(iv) If \ \ $x_i \prec x_j \prec x_k\ \ $ and the pair \ \ $(c_{ij},
d_{ij})\ \ $ is replaced by a
scalar multiple so that in (iii) the scalar \ \ $\lambda\ \ $ is unity, then,
for some
scalar \ \ $\mu \in \IC_*\ \ $ we have \ \ $c_{ik} = \mu c_{jk}\ \ $ and \ \
$d_{ik} = \mu d_{jk}.\ \ $

Suppose further, that  \ \ $X \subseteq Y\ \ $ is a subset such that the
natural map \ \ $H_1 (\Delta (X))
\rightarrow H_1 (\Delta (Y))\ \ $ is the zero map. Then there exists a choice
of
elements \ \ $g_i\ \ $ in \ \ $G_i,\ \ $ for all \ \ $x_i \in X,\ \ $ such that
for each pair we have \ \ $(c_{ij},
d_{ij}) = (\lambda_{ij} g_i \lambda_{ij} g_j)\ \ $ for some \ \ $\lambda_{ij}
\in \IC_*.\ \ $
\vspace{.3in}\rm

\begin{mproof} Let \ \ $(X, \prec)\ \ $ be viewed as a graph. We may assume
that it is
connected. Let \ \ $\tau\ \ $ be a maximal tree in \ \ $X.\ \ $ Fix an edge \ \
$(x_1, x_k)\ \ $ in
\ \ $\tau\ \ $ and define \ \ $g_1 = d_{1, k}.\ \ $ Using the
edges of \ \ $\tau,\ \ $ and properties
(i), (ii) and (iii), define \ \ $g_i\ \ $ recursively  for all the vertices \ \
$x_i\ \ $
of \ \ $X.\ \ $ In the process of doing this, whenever property (iii) is used
replace
the `new' pair by a scalar multiple so that the scalar \ \ $\lambda\ \ $ of
(iii) is
unity.

Let \ \ $C_1 (\Delta (Y))\ \ $ be the group of 1-chains of the complex \ \
$\Delta (Y).\ \ $ If
\ \ $(x_i, x_j)\ \ $ (with \ \ $x_i \prec x_j\ \ $) is a 1-simplex of \ \
$\Delta (Y)\ \ $ then by (i), (ii)
and (iii) we know that

\[
(c_{ij}, d_{ij}) = (\alpha g_i, \beta g_j)
\]

\noindent for some scalars \ \ $\alpha, \beta \in \IC_*.\ \ $ Define the group
homomorphism \ \ $\Phi:
C_1 (\Delta (Y)) \rightarrow \IC_*\ \ $ by taking the unique extension of the
correspondences \ \ $\Phi ((x_i, x_j)) = \alpha \beta^{-1}.\ \ $ By (iv) we
have \ \ $\Phi
(\partial \sigma) = 1\ \ $ if \ \ $\sigma\ \ $ is a 2-simplex of \ \ $\Delta
(Y).\ \ $ It
follows then that \ \ $\Phi (w) = 1\ \ $ whenever \ \ $w\ \ $ is a 1-boundary.

Let \ \ $(x_i, x_j)\ \ $ be an edge of \ \ $(X, \prec)\ \ $ which is not an
edge of \ \ $\tau.\ \ $ We
are required to prove that \ \ $\alpha = \beta\ \ $ in this case. But since \ \
$\tau\ \ $ is a
maximal tree, there is a 1-cycle \ \ $w\ \ $ for \ \ $\Delta (X)\ \ $
consisting of
the simplex \ \ $(x_i,x_j)\ \ $ and (distinct) 1-simplexes from \ \ $\tau.\ \ $
Thus \ \ $\Phi (w) = \Phi ((x_i, x_j))
= \alpha \beta^{-1}.\ \ $ By the hypothesis \ \ $w\ \ $
is a boundary, and so, using the previous
paragraph, \ \ $\alpha = \beta.\ \ $
\end{mproof}

\noindent {\bf The proof of Theorem 4.2}

Let \ $\ \cal L\ $\  be represented as in Theorem 2.2. Let \ $\ \cq
= \{U_1, \dots,
U_n\}\ $\  be a partition of \ $\ M({\cal L})\ $\  generated by basic clopen
sets,
and choose \ $\ x_i \in U_i\ $\  for each \ $\ i.\ $\  Let \ $\ Q_i =
\tilde{P}(U_i)\ $\  and
associate with
\ $\ \cq\ $\  the subalgebra \ $\ {\cal
A} (\cq)\ $\  of \ $\ \cal A\ $\  given by

\[
{\cal A} (\cq) = \hbox{ span } \{Q_i {\cal L} (\ch)
Q_j : Q_i \ll Q_j \}
\]

\noindent where \ $\ Q_i \ll Q_j\ $\  if and only if there exists \ $\ Q_k\ $\
with \ $\ Q_i \prec Q_k \prec
Q_j\ $\ . By Lemma 4.3 for each such pair there exists an operator pair
\ $\ (C_{ij}, D_{ij})\ $\  with \ $\ C_{ij} \in Q_i \cc, D_{ij} \in Q_j \cc\ $\
, such
that \ $\ \alpha (A) = C_{ij} A D_{ij}^{-1}\ $\  for all operators \ $\ A\ $\
with \ $\ A =
Q_i A Q_j\ $\ .

We show that the restriction \ $\ \alpha | {\cal A}(\cq)\ $\  is inner. Let \
$\ (X,
\ll) \ $\  be the set \ $\ X\ $\  with the partial order inherited from \ $\
\cq\ $\ . View
the associated complex \ $\ \gD(X, \ll)\ $\  as a subcomplex of the
topological complex \ $\ \gD_{top}(M(\cl), \le)\ $\ .
Then, by the spectral homology hypothesis the natural inclusion induced
map

\[
H_1 (\gD(X, \ll)) \rightarrow H_1^{sp}(\ca)
\]

\noindent is the zero map and so there exists a finite subset \ $\ Y_0\ $\  of
\ $\ M({\cal
L})\ $\  containing \ $\ X\ $\  such that the natural map

\[
H_1 (\gD(X, \ll)) \rightarrow H_1 (\gD_{top}(Y_0, \le))
\]

\noindent is the zero map. By Theorem 2.5 there is a finer partition associated
with a finite set \ $\ Y\ $\  containing \ $\ Y_0\ $\  such that
$\ \gD_{top}(Y_0, \le)$\ is a subcomplex of $\ (\gD(Y, \ll)$.\
Thus, taking compositions, it follows that the
natural map

\[
H_1 (\gD(X, \ll)) \rightarrow H_1 (\gD(Y, \ll))
\]
is the zero map.

Recall, from the discussion preceeding Lemma 4.5 that the pairs of
operators \ $\ (C_{ij}, D_{ij})\ $\  that are associated with the partition
satisfy the requirements (i), (ii), (iii), (iv) of Lemma 4.5. It follows
from the lemma that there is a choice of non-zero scalars \ $\ \lambda_{ij}\ $\
and invertible operators \ $\ C_i\ $\  in the algebra \ $\ Q_i \cc\ $,
\  for \ $\ i = 1,
\dots, m,\ $\  such that

\[
(C_{ij}, D_{ij}) = (\lambda_{ij} C_i, \lambda_{ij} C_j)
\]

\noindent for all \ $\ Q_i \ll Q_j\ $\ . It now follows that \ $\ \alpha |
{\cal A} (\cq)\ $\  is
inner.

Finally we show that \ $\ \alpha\ $\  is a point weak star limit of inner
automorphisms.

Consider a chain of partitions $\cq_1 \subseteq \cq_2 \subseteq \dots$
with associated subalgebra chain \ $\ \{{\cal A} (\cq_k)\}.\ $\
Then the restriction of $\ \alpha\ $\  to $\ {\cal A} (\cq_k)$\ \ has
the form \\ $\ A \to C_kAC_k^{-1}$,\  for some invertible operator $\ C_k$\
in $\cc$.
Let

\[
P_k = \mbox{ sup}\{Q_i, Q_j : Q_i \ll Q_j, \ \ Q_i, Q_j \in \cq_k\}.
\]

\noindent Then we can arrange that $\ C_{k+1}$\ extends $\ \ C_{k}$\  \
in the sense that $\ \ C_{k+1}P_k = C_kP_k.$\ \
To be more precise about this, let us restrict attention in the
remainder of the proof to the case of $\ \ca$\  irreducible. The general
case then follows readily. By Proposition 5.1 of
\cite{scp-tensor} $\ (X, \le)$\  is a connected binary relation.
It follows that having chosen $\ C_1$\ subsequent choices of the
$\ C_k$\  are naturally uniquely determined on the projection
$\ P_k$.\  One way to see this is to note
that Theorem 2.6 implies that
projections  $\ Q_i \in \cq_k$ \ and $\ Q_j \in \cq_1$ \
have subprojections
$E\ $ and $F,\ $ respectively, which are
$\ll$ -connected.
Note that it follows that the bounded operator $\ C_k$\ implements
\ $\alpha$\  on the space $\ P_k\ca P_k$\  as well as on the smaller
subalgebra $\ \ca(\cq_k)$.\

Suppose for the moment that the union of the algebras
 \ $\ {\cal A} (\cq_k)\ $\ is weak star dense in $\ca$. We show that
if $\ A \in \ca$\ and $\  \phi \ $ is a weak star continuous
functional, then for $\ \epsilon > 0$\ there exists an invertible
operator $\ D$\ in $\ \cc$\ such that

\[
\ |\phi(\alpha(A)) - \phi(DAD^{-1})\ | < \epsilon.
\]

\noindent In view of the hypothesised weak star density,
the increasing projections $\ P_k$\  converge to the identity operator
in the strong operator topology.
It follows that there is a densely defined (possibly) unbounded
operator $\ C$\ whose domain is the linear span of the ranges of the
projections \ $P_k$,\
such that $\ CP_k = C_kP_k$.\  Let $\ E_k^1, E_k^2 $\ and \ $E_k^3$\ \ be the
spectral projections
for $\ |C|$\ for the sets
$\ [0,k^{-1}), \ [k, \infty)$\ \ and $\ [k^{-1},k)$\ \ respectively.
These are projections in $\ \cc$\  \ and the sequences
$\ E_k^1$\ \ and $\ E_k^2$\  \ converge to zero in the strong operator
topology. Consider the sequence of operators $\ D_k$\ in $\ \cc$\ given by

\[
D_k = k^{-1}E_k^1 + kE_k^2 + CE_k^3.
\]

\noindent Let $\ \phi$\ be a weak star continuous linear linear functional
on $\ \ca$, \ so that $\ \phi(A) = $\ trace$(TA)$\ \
for some trace class operator $\ T$.\  For
fixed $\ A \in \ca$ \ and $\ \epsilon > 0$\ \
choose $\ k_0$\ large enough so that

\[
\ |\phi(\alpha(A)) - \phi(\alpha(E_k^3AE_k^3))\ | < \epsilon/2
\]

\noindent for all $\ k > k_0$. (Recall that $\alpha$ is automatically weak star
continuous.) We shall show that by increasing $k$, if necessary, we have

\[
\ |\phi(D_kAD_k^{-1}) - \phi(D_kE_k^3AE_k^3D_k^{-1})\ | < \epsilon/2,
\]

\noindent from which the desired conclusion follows, since
\ $D_kE_k^3AE_k^3D_k^{-1} = \alpha(E_k^3AE_k^3)$.\ \
It suffices to show that for large $k$

\[
\ |\mbox{trace}(XD_kE_k^iAE_k^jD_k^{-1}\ | < \epsilon/16
\]

\noindent for the eight pairs $\ (i,j)$ with $(i,j) \ne (3,3)$.\
For $\ (i,j) = (1,1)$\ or $\ (2,2)$\ the quantity is simply
$\ \ |$trace$(XE_k^iAE_k^j)\ |
$\ \
and so these cases are clear. For
the case $\ i = 2$\ and $j = 3$\ observe that

\[
D_kE_k^2AE_k^3D_k^{-1} = kE_k^2AC^{-1}E_k^3 =
(C_1CE_k^2)E_k^2AE_k^3(C^{-1}E_k^3) = C_1\alpha(E_k^2AE_k^3)
\]

\noindent where $\ C_1$\ is a contraction such that $\ C_1CE_k^2 = kE_k^2$.\
Consequently
this case is also clear,
since $\ XE_k^2 \to 0 $\  in the trace class norm
and the equality above shows that the operators
$\  D_kE_k^2AE_k^3D_k^{-1}$\ are uniformly bounded by $ \|\alpha\| \|A\|.$\
The case $\ (i,j) = (1,3)$\  is similar. The cases $\ (3,2), (3,1) $\ and
$\  (1,2)\  $ are elementary,
and so it remains to consider
$\ i = 2 $ and $ j = 1.$\
In this case it can be seen that for $\ k > \|\alpha\|^2$\  the operator
$\ E_k^2AE_k^1$\ is zero.
This need only be observed for rank one operators $\ A$ in $\ca,$\  \
since their linear span is
dense in the weak operator topology, and, in view of Lemma 4.4 (i),
this verification is elementary.

To finish the proof we confirm the technical detail that a subalgebra
chain $\ \ca(\cq_k)$\  can be found, with dense union.

In view of the separability of the underlying Hilbert space
there is a countable family $\ L_1, L_2, \dots$\  such that each projection
$\ L$\ in $\ \cl\ $ is both the supremum and infimum of projections in the
family $\ \{L_k\}$.\
Let $\ \cq_k$\ \ be the partition of $\ X$\  generated by $\ L_1,\dots L_k$.\
It will be enough to show that each rank one operator
$\  R $\ \ is in the norm closure of
the associated algebras $\ \ca(\cq_k)$.\

By Lemma
4.4, $\ R = x \ot y^*$\ and there is a projection \ $\ L\ $\  in \ $\ \cal L\
$\  such
that \ $\ Lx  = x\ $\  and
\ $\ (I-L_-)y =y\ $.
By norm approximation we may reduce to the case
that both $\ L$\  and $\ L_-$\
belong to $\ \{L_k\}$.\
We find a particular projection \ $\ L_\epsilon \leq L\ $,\
with $\ L_\epsilon $ also in \ $\ \{L_k\}$,\
so that
if \ $\ x_\epsilon = L_\epsilon x\ $\  then
\ $\ R_\epsilon = x_\epsilon \otimes y^*\ $\  is a rank one
operator close to \ $\ R\ $\  and lying in one of the algebras
$\ \ca(\cq_k)$.\
We do with the following
argument borrowed from the proof of
Lemma 2.7 of \cite{orr-scp}.

Since \ $\ \cal L\ $\  is completely distributive we have, by
Lemma 2.3 of \cite{orr-scp} for example,

\[
L = \sup \{G_+ : L \not\le G_{+-}, G \in {\cal L}\}
\]

\noindent for all \ $\ L\ $\  in  \ $\ \cal L\ $\  with \ $\ L \neq 0.\ $\
Readjusting our choice of $\ \{L_k\}$\ if necessary,
we may assume that $\ L$\  is in fact the supremum of
projections $\ G_+$\  which belong to $\ \{L_k\}.$\
Thus we may choose  \ $\ G^1_+, \dots
G^n_+\ $\  in $\ \{L_k\}$\
so that the projection \ $\ L_\epsilon =
$\ sup$\{G^i_+ : 1 \le i \le n\}$
determines the rank one operator \ $\ R_\epsilon\ $\
as above with \ $\ ||R_\epsilon - R|| < \epsilon \ $\ . Note
that \ $\ [H, I)^c \prec [H_-, I)\ $\  for
any projection \ $\ H\ $\  and so

\[
[G^i_+, I)^c \ \prec \ [L, I)^c \cap [G^i_{+-}, I)\  \prec \ [L_-, I)
\]

\noindent for each \ $\ i\ $. Furthermore the middle sets here are nonempty.
Thus by writing $\ L_\epsilon x$\ as a sum of vectors in
the ranges of the projections $\ G^i_+$\  we see that
\ $\ R_\epsilon\ $\  belongs to the algebra \ $\ {\cal A} (\cq)\ $\
associated
with the partition generated by \ $\ [L_-, I)\ $\
and \ $\ [G^i_+, I)^c, 1 \le i
\le  n\ $\ . \hfill \ $\ \Box\ $\
\vspace{.3in}

\noindent {\bf The proof of Theorem 4.1}

\noindent Let \ $\gd\ $ be a weak star continuous derivation of \ $\ca\ $.
We wish to show that
\ $\gd \in B^1_{ess}(\ca)\ $. By a standard argument of Kadison
and Ringrose \cite{kad-rin} we may
assume that \ $\gd(C) = 0\ $ if \ $C \in \cc\ $ where \ $\cc\ $ is a fixed masa
in \ $\ca\ $. Let \ $\ga\ $ be the automorphism \ \ exp\ $(\gd)\ $
of \ $\ca\ $. Calculation shows
that it is a Schur automorphism of \ $\ca\ $ with respect to \ $\cc\ $.

Let \ $\ca(\cq) \subseteq \ca\ $ be a subalgebra associated
with a finite partition as in
the last proof. Then, by the arguments above,
the restriction \ $\ga_{\cq} = \ga |\ca(\cq)\ $ is inner and has
the form \ $\ga_C\ $ where \ $\ga_C(A) = CAC^{-1}\ $
for all \ $A\ $ in \ $\ca(\cq)\ $, where \ $C\ $ is an invertible element of \
$\cc\ $.
Let \ $D\ $ be a logarithm for \ $C\ $ in \ $\cc\ $ associated with the inner
derivation
\ $\gd_D\ $. Then \ \ exp\ $(\gd_D) = \ga_{e^D} = \ga_C = \ $exp\ $(\gd_{\cq})\
$,
where \ $\gd_\cq\ $ is the restriction \ $\gd |\ca(\cq)\ $. (Our assumption for
\ $\gd\ $
implies that the weakly closed \ $\cc\ $-bimodules in\ $\ca\ $ are
invariant for \ $\gd\ $.) Since \ $\gd_\cq\ $ commutes with \ $\gd_D\ $ we
deduce that
\ \ exp\ $(\gd_D - \gd_\cq) =\ $ id\ \ , and hence that \ $\gd_D = \gd_\cq\ $.
It now follows, as in the last proof, by the density of
the algebras \ $\ca(\cq)\ $, that \ $\gd\ $ is a point
weak star inner automorphism.
\vspace{.3in}

\noindent {\bf Remark 4.6}\ \  We conjecture that the higher order analogue of
Theorem 4.1 also holds, that is,
\[
H^{sp}_n ({\cal A}) = 0 \ \Rightarrow \hbox{ Hoch}_{ess}^n ({\cal A}) =0.
\]
This is known to be true in the case of digraph algebras.
In fact if \ $H_n^{sp}(A(G)) = 0 \ $ then \ $H_n(\gD(G_r)) = 0$,
by Proposition 3.7. Thus \ $H^n(\gD(G_r)) = 0\ $ by the duality for simplicial
cohomology, and so Hoch$^n(A(G)) = 0\ $ by the cohomological identifications
given in Gerstenhaber and Schack \cite{ger-sch} and Kraus and Schack
\cite{kra-sch}.

The analysis above  complements some of the results of
Gilfeather and Moore in \cite{gil-moo}. They have shown, in particular,
that if \ \ $\ga\ \ $ is an automorphism of
{\it any} CDCSL algebra \ \ $\ca\ $ then the following conditions are
equivalent.

(i) rank($\ga(R)) = \ $rank$ (R)\ $ for all finite rank operators in \ $\ca.\ $

(ii) \ $\ga\ $ is {\it quasispatial} in the sense that there is a
 closed injective linear transformation
\ $T : \ch_1 \to \ch_2\ $, whose range and domain are dense,
such that  \ $\ga(A)Ty = TAy\ $ for all \ $y\ $ in the domain of \ $T.\ $

Let \ $\ga\ $ be a Schur automorphism, which is pointwise weakly inner in
the sense of Theorem 4.2. Then \ $\ga\ $ must preserve the rank of finite rank
operators
and so it follows from Theorem 4.2 and the Gilfeather-Moore result that if
\ $H^{sp}_1(\ca) = 0 \ $ then all Schur automorphisms are quasispatial.
In fact we can see this directly in the proof of Theorem 4.2
where the possibly unbounded implementing operator
for the automorphism is constructed.

In the other direction it seems plausible that if \ $H^{sp}_1(\ca) \ne 0 \ $
then there exists a Schur automorphism of \ $\ca\ $ which is not quasispatial.
More generally, it would be interesting to determine whether there is a
converse implication
to the one conjectured above.

\section{Final Remarks}

  It would be interesting to see to what extent
it is possible to develop a homology theory for general CSL
algebras which is
based upon spectral invariants. One difficulty apppears to be that the
unitary invariants for general commutative projection
lattices, as given by Arveson in \cite{arv}
for example, do not have the explicit intrinsic nature as those
of Theorem  2.2. Ideally one would wish to develop a general
theory capable of calculations of the homology of even infinite
tensor products.
Another natural direction is, of course, to extend the spectral
homology invariants to
other classes of reflexive operator algebras. One can envisage that
such a development is possible for classes of projection
lattices generated by
CDCSL lattices together with "amenable" lattices, which are not necessarily
commuting,
through the natural operations -  products, fibre sums, joins, and so forth
- at the algebraic level.
For example, if $\ \ \ca_1$\ \  is a CDCSL algebra and \ \ $ca_2\ \ $
is the reflexive
algebra associated with two projections in generic position,
(see Lambrou and Longstaff \cite{lam-lon} for example) then
the
vanishing of Hochschild cohomology
of \ \ $\ \ \ca_1 \ot \ca_2\ \ $\ \
should be a consequence of vanishing spectral homology.

In \cite{scp-stable} we develop a stable homology
theory for general \nsa \ operator algebras which is based
on partial isometries normalising a given masa.
Roughly speaking, if \ $\ \ca\ $\  is an operator algebra with masa
\ $\ \cc\ $\  then the stable homology
group \ $\ H_1(\ca;\cc)\ $\  is an abelian
group associated with cycles of partial isometries in the algebras
\ $\ M_n \ot \ca\ $\  which
normalise \ $\ D_n \ot \cc.$  By hypothesis such partial
isometries must form part of a complete
matrix unit system.
The precise definition of  \ $\ H_n(\ca;\cc)\ $\  is analogous in spirit to the
definition of
the \ $\ K_0\ $\ -group. In fact  if the inclusion
\ $\ \cc \to C^*(\ca)\ $\  induces a regular surjection
\ $\ K_0\cc \to K_0(C^*(\ca))\ $\  then \ $\ H_0(\ca;\cc)$\   is often
equal to $\ K_0(C^*(\ca))$.

In a CDCSL algebra a masa is unique up to inner unitary equivalence,
and so in this case   stable homology is an invariant for the algebra
and we may denote the groups
simply as  \ $\ H_n(\ca$).\ \
The algebra \ $\ {\ell^\infty}\ $\  of diagonal
operators relative to an orthonormal basis is a CDCSL algebra
and in this case
the first stable homology group  \ $\ H_0(\ell^\infty) $ is equal to
\ $\  K_0(\ell^\infty).\ $\
On the other hand \ $\ H_n^{sp}(\ell^\infty)\ $\  is
simply the restricted direct product \ $\ \ZZ^\infty.\ $\

Another contrast, of a different and perhaps more significant nature,
is that \ $\ H_*^{sp}(\ca) \ $\
is computed purely in terms of the structure of
\ $\ \cl = \ $\ Lat$\ca,\ \ $\ \ and so spectral homology
takes no acccount of the
fact that \ \ $ cl$\ \  may have atoms of both
finite and infinite rank. Stable homology on the
other hand is based on cycles of partial isometries of the same rank
(in the case of CSL algebras) and so from this point of view
provides a more discriminating invariant.

\end{document}